\documentclass[12pt]{article}
\usepackage{amssymb,amsmath}

\newtheorem{Theorem}{Theorem}
\newtheorem{Example}{Example}
\newtheorem{Definition}{Definition}
\newtheorem{Lemma}{Lemma}
\newenvironment{proof}{
  \noindent\textbf{Proof}\ }{\hspace*{\fill}
  \begin{math}\Box\end{math}\medskip}

\begin{document}

\title{On chromatic number of unit-quadrance graphs (finite Euclidean graphs)}
\author{Le Anh Vinh\\
School of Mathematics\\
University of New South Wales\\
Sydney 2052 NSW} 
\date{\empty}
\maketitle

\begin{abstract}
  The quadrance between two points $A_1 = ( x_1, y_1 )$ and $A_2 = ( x_2, y_2
  )$ is the number $Q ( A_1, A_2 ) = ( x_1 - x_2 )^2 + ( y_1 - y_2 )^2$. Let
  $q$ be an odd prime power and $F_q$ be the finite field with $q$ elements.
  The unit-quadrance graph $D_q$ has the vertex set $F_q^2$, and $X, Y \in
  F_q^2$ are adjacent if and only if $Q ( A_1, A_2 ) = 1$. Let $\chi ( F_q^2 )$ be
  the chromatic number of graph $D_q$. In this note, we will show that $q^{1 /
  2} ( 1 / 2 + o ( 1 ) ) \leqslant \chi ( F_q^2 ) \leqslant q ( 1 / 2 + o ( 1
  ) )$. As a corollary, we have a construction of triangle-free graphs $D_q$ of order $q^2$ with $\chi ( D_q ) \geq q/2$ for infinitely many values of $q$.
\end{abstract}

\section{Introduction}

In \cite{quadrance}, Wildberger introduces a remarkable new approach to trigonometry and Euclidean geometry by replace distance by quadrance and angle by spread, thus allowing the development of Euclidean geometry over any field. The following definition follows from \cite{quadrance}.

\begin{Definition}
  The \textbf{quadrance} $Q ( A_1, A_2 )$ between the points $A_1 = ( x_1,
  y_1)$, and $A_2 = ( x_2, y_2 )$ in $F_q^2$ is the number
  \[ Q ( A_1, A_2 ) := ( x_2 - x_1 )^2 + ( y_2 - y_1 )^2 . \]
\end{Definition}

This approach motives the question of finding chromatic number of unit-quadrance graph over any field. Note that in usual $2$-dimensional Euclidean space $R^2$ then the quadrance between $A_1, A_2$ is unit if and only if the distance between $A_1, A_2$ is unit. 

Let $q$ be an odd prime power and $F_q$ be the finite field with $q$ elements. The unit-quadrance graph $D_q$ has the vertex set $F_q^2$, and $X, Y \in F_q^2$ are adjacent if and only if $Q ( A_1, A_2 ) = 1$. This graph (or so-called finite Euclidean graph) was also studied by Medrano et al in \cite{before}. Recall that the graph $D_q$ is a regular graph with degree $\Delta(D_q) = q - (-1)^{(q-1)/2}$ (see \cite{before}).  Let $\chi ( F_q^2 )$ be the chromatic number of graph $D_q$. The main result of this note is the following theorem.

\begin{Theorem}\label{main}
  Suppose that $q = p^n > 3$ where $p$ is an odd prime number then
  \[ q^{1 / 2} ( 1 / 2 + o ( 1 ) ) \leqslant \chi ( F_q^2 ) \leqslant
     \frac{p^n + p^{n - 1}}{2} = q ( 1 / 2 + o ( 1 ) ) . \]
\end{Theorem}

\section{Some Lemmas}

In order to prove Theoreom 1, we need some lemmas.

\begin{Lemma}\label{line}
Suppose that $a \in F_q$ such that  $a^2 + 1$ is not square in $F_q$. For any $A\neq B$ in the line $y=ax + i$ then $Q(A,B) \ne 1$.
\end{Lemma}

\begin{proof}
  Suppose that $A = ( x_1, a x_1 + i )$ and $B = ( x_2, a x_2 + i )$ for some $x_1 \neq x_2
  \in F_q$. We have
  \[ Q ( A, B ) = ( a^2 + 1 ) ( x_1 - x_2 )^2 \neq 1 \]
  since $a^2 + 1$ is not square in $F_q$. The lemma follows.
\end{proof}

Recall that a (multiplicative) character of $F_q$ is a homomorphism from
$F_q^{\ast}$, the multiplicative group of the non-zero elements of $F_q$, to the
multiplicative group of complex numbers with modulus $1$. The identically $1$
function is the principal character of $F_q$ and is denoted $\chi_0$. Since
$x^{q - 1} = 1$ for every $x \in F_q^{\ast}$ we have $\chi^{q - 1} = \chi_0$
for every character $\chi$. A character $\chi$ is of order $d$ if $\chi^d =
\chi_0$ and $d$ is the smallest positive integer with this property. By
convention, we extend a character $\chi$ to the whole of $F_q$ by putting
$\chi ( 0 ) = 0$. The quadratic (residue) character is defined by $\chi ( x )
= x^{( q - 1 ) / 2}$. Equivalently, $\chi$ is $1$ on square, $0$ at $0$ and $-
1$ otherwise. It is easy to see that $\sum_{i \in F_q} \chi ( i ) = 0$. We have the following lemma.

\begin{Lemma}\label{exist} Suppose that $t$ is not a square in $F_q$. Then there exists a square number $i \in F_q^\ast$ such that $-i+t$ is not a square in $F_q$.
\end{Lemma}

\begin{proof} Let $A_q$ be the number of $i \in F_q^\ast$ such that $i$ is square but $-i+t$ is not.
Consider the sum
  \begin{align*}
    \sum_{i \in F_q} \left( 1 - \chi(i) \right) \left( 1
    + \chi(-i+t) \right) &= q + \sum_{i \in F_q} \left\{ - \chi(i) + \chi(-i+t) 
    - \chi(i ( -i + t )) \right\}\\
    &= q - \sum_{i \in F_q \neq 0, t} \chi(i^2 ( -1 + ti^{- 1} ))\\
    &= q - \sum_{i \in F_q \neq 0, t} \chi(-1 + ti^{- 1})\\
    &= q + \chi(-1) + \chi(0)\\   
    &= q + (-1)^{(q-1)/2}. 
  \end{align*}
But we have
  \begin{equation*}
   \left( 1 - \chi(i) \right) \left( 1 + \chi(-i + t)\right) =
   \begin{cases}
    4  & \text{if}\; \; i \;  \text{is square} \;\text{and}\; -i+t \; \text{is not square}  ,\\
    2  & \text{if}\; \; i = t.\\
    0 & \text{otherwise}.
   \end{cases}
  \end{equation*}
  Thus, $4 A_q + 2 = q + (-1)^{(q-1)/2}$ or $A_q = ( q + (-1)^{(q-1)/2} - 2 ) / 4\geq 1$. This concludes the proof of the lemma.
\end{proof}  

\begin{Lemma}\label{two lines} Suppose that $q \geq 5$ is a prime power. Let $a \in F_q$ such that $a^2 + 1$ is not square. There exist $t \in
  F_q^{\ast}$ such that if $A$ is in the line $y = a x + i$ and $B$ is in the
  line $y = a x + i + t$ then $Q ( A, B ) \neq 1$ for any $i \in F_q$.
\end{Lemma}

\begin{proof}
  Suppose that $X = ( x, a x + i )$ and $Y = ( y, a y + i + t )$ then
  \begin{align*}
    Q ( X, Y ) &= ( x - y )^2 + ( a ( x - y ) + t )^2\\
    &= ( a^2 + 1 ) ( x - y )^2 + 2 a ( x - y ) t + t^2 .
  \end{align*}
  If $Q ( X, Y ) = 1$ then $a^2 + 1 = [ ( a^2 + 1 ) ( x - y ) + a t ]^2 +
  t^2$. From Lemma \ref{exist}, there exists $i \in F_q^\ast$ such that $i$ is square but $-i+a^2+1$ is not. Hence, we can choose $t$ such that $t^2 = i$. This concludes
  the proof of the theorem.
\end{proof}

In \cite{before}, Medrano et al. give a general bound for eigvenvalue of $D_q$.

\begin{Lemma}\label{eigenvalue} (\cite{before}) Let $\lambda \ne \Delta(D_q)$ be any eigenvalue of graph $D_q$ then $| \lambda | \leq q^{1/2}$.
\end{Lemma}

The following result which is due to Hoffman gives us a connection between the chromatic number and eigenvalues of a graph.

\begin{Lemma}\label{bound} (\cite{hof}) Let $G$ be any graph with the largest and least eigenvalues are $\lambda_1$ and $\lambda_{\ast}$. Then \[\chi(G) \geq 1 - \frac{\lambda_1}{\lambda_{\ast}}.\]
\end{Lemma}

\section{Proof of Theorem 1}

The lower bound is straigtforward from Lemmas \ref{eigenvalue} and \ref{bound}. We have $\lambda_1 = q \pm 1$ since $D_q$ is a $(q\pm 1)$-regular graph. Lemma \ref{eigenvalue} gives us $\lambda_{\ast} \geq -q^{1/2}$. Hence 
\[\chi(D_q) \geq 1 + \frac{q \pm 1}{q^{1/2}} = q^{1/2}(1+o(1)).\]

The upper bound is proved by a geometrical colouring. Let $a \in F_q^{\ast}$ such that $a^2 + 1$ is not a square in $F_q$ and choose $t$ satisfy Lemma \ref{two lines}. We partition $F_q$ into $p^{n-1}$ disjoint sets of the form $\{\alpha_i, \alpha_i +t,\ldots, \alpha_i +(p-1)t\}$ for $i = 1,\ldots, p^{n-1}$. We colour $2p$ points in the lines $y = a x + \alpha_i + (2k)t$ and $y = a x + \alpha_i + (2k+1)t$ the same colour for each $i = 1, \ldots, p^{n-1}$ and $0 \leq k \leq (p-3)/2$ (using different colours for each pair $(i,k)$). And we colour $p$ points in the line $y= a x + \alpha_i + p-1$ the same colour for each $i=1,\ldots, p^{n-1}$ (using different colour for each $i$). Then total number of colours we have used is $p^{n-1}(p+1)/2$. From Lemmas \ref{line} and \ref{two lines}, there does not exist two point $A, B$ such that $Q(A,B)=1$ and $A, B$ have the same colour. This implies the upper bound and concludes the proof of the theorem.

\begin{Example}
  For $q = 7$ then we have $3 \leq \chi(F_7^2) \leq 4$ (the upper bound follows from Theorem \ref{main}, the lower bound is trivial since $D_7$ contains a cycle of length $7$). We choose $a = 5, t = 3$ then $a^2+1 = 5$ and $-t^2 + a^2 + 1 = 3$ are not square in $F_7$. We have a $4$-colouring of points of $F_7^2$ in Table 1. It has been verified by computer that there does not exist a $3$-colouring of vertices of $D_7$. This implies that $ \chi(F_7^2) =  \chi(D_7) = 4$. 
\begin{table}
	\centering
 \begin{tabular}{|c|c|c|c|c|c|c|}
    \hline
    3&1&2&3&4&1&2\\
    \hline
    2&3&4&1&2&3&1\\
    \hline
    4&1&2&3&1&2&3\\
    \hline
    2&3&1&2&3&4&1\\
    \hline
    1&2&3&4&1&2&3\\
    \hline
    3&4&1&2&3&1&2\\
    \hline
    1&2&3&1&2&3&4\\
    \hline
  \end{tabular} 
	\caption{$4$-colouring of $F_7^2$}
	\label{q7}
\end{table}
  
\end{Example}

\section{Triangle-free graphs with arbitrary high chromatic numbers}

The unit-quadrance graphs give us a family of triangle-free graphs with arbitrary high chromatic numbers. We have the following lemma.

\begin{Lemma}\label{no triangle}
Let $q$ be any prime of the form $q = 12k \pm 7$. Then $D_q$ contains no triangle.
\end{Lemma}

\begin{proof}
Suppose that $D_q$ contains a triangle $XYZ$ with $X = ( m, n )$, $Y = ( m + x, n + y )$ and $Z = ( m + x + u, n + y + v )$  for some $m, n, x, y, u, v \in F_q$ then $x^2 + y^2 = u^2 + v^2 = (x+u)^2 + (y+v)^2 = 1$. This implies that $x u + y v = -1/2$.
We have $( x u + y v )^2 + ( x v - y u )^2 = ( x^2 + y^2 ) ( u^2 + v^2)$ so
  \[ ( x v - y u )^2 = 1 - 1/4 = 3/4. \]
For $q = 12 \pm 7$ then $3$ is not square in $F_q$. Hence we have a contradiction. This concludes the proof of the lemma.
\end{proof}

From Theorem \ref{main} and Lemma \ref{no triangle}, if $q \equiv \pm 7$ (mod $12$) is a prime number then $D_q$ is a triangle-free graph with chromatic number $\chi(D_q) \geq q/2(1+o(1))$. This give us a family of triangle-free graphs with arbitrary high chromatic numbers.

\section{Higher dimensional spaces}

The quadrance can be defined in higher dimension space $F_q^m$ for $m \geq 2$ as follows. 

\begin{Definition}
  The \textbf{quadrance} $Q ( A_1, A_2 )$ between the points $X = (x_1,\ldots,x_m)$, $Y = ( y_1,\ldots,y_m )\in F_q^m $ is the number
  \[ Q ( A_1, A_2 ) := \sum_{i=1}^{m} ( x_i - y_i )^2. \]
\end{Definition}

The unit-quadrance graph $D_q^{m}$ has the vertex set $F_q^m$, and $X, Y \in
F_q^m$ are adjacent if and only if $Q ( A_1, A_2 ) =1$. Let $\chi ( F_q^m )$ be
the chromatic number of graph $D_q^m$. Similar as the above, we have the following theorem.

\begin{Theorem}\label{gen}
  Suppose that $m\geq 2$ and $q = p^n > 3$ where $p$ is an odd prime number then
  \[ q^{(m-1) / 2} ( 1 / 2 + o ( 1 ) ) \leqslant \chi ( F_q^m ) \leqslant
     \frac{q^{m-2}[p^n + p^{n - 1}]}{2} = q^{m-1} ( 1 / 2 + o ( 1 ) ) . \]
\end{Theorem}

The proof of this theorem is omitted since it is the same as the proof of Theorem 1.

\end{document}